\documentclass[11pt]{article} 
  \usepackage{epsfig} 
  \usepackage{amsmath} 
  \usepackage{amssymb} 
  \usepackage{amsthm} 
  \usepackage{graphicx}
  \usepackage{psfrag}
\newcommand\inta{\int^{\infty}_{-\infty}}
\newcommand\into{\int^{\infty}_{0}}

\newcommand\s{Schr\"odinger }

\renewcommand\Im{ {\mathrm{Im}\!\!\!} \ \  }

\begin{document}  
\title{Boundary Value Problems for Linear PDEs with \\
Variable Coefficients}
\author{A.S. Fokas \\
{\em Department of Applied Mathematics and Theoretical Physics} \\
{\em University of Cambridge} \\
{\em Cambridge CB3 0WA, UK}}
\date{June 2003}

\maketitle

\centerline{\bf Abstract}

\vskip .2in

A new method is introduced for studying boundary value problems for a
class of linear PDEs with {\it variable}
coefficients.  This method is based on ideas recently introduced by the author
for the study of boundary value problems for PDEs with {\it constant}
coefficients.  As illustrative examples the following boundary value problems
are solved: (a) A Dirichlet and a Neumann problem on the half line for
the time-dependent \s equation with a space dependent potential.  (b)
A Poincar\'e problem on the quarter plane for a variable coefficient
generalisation of the Laplace equation.

\pagebreak

\section{Introduction}
The main aim of this paper is to introduce a methodology for solving boundary
value problems for a class of linear PDEs with variable coefficients.  As
illustrative examples the following boundary value problems will be solved:

(i) A Dirichlet problem for the time-dependent \s equation on the half
line,
$$ iq_t + q_{xx} + u(x)q=0, \quad 0<x<\infty, \quad t>0, \eqno (1.1)$$
$$ q(x,0) = q_0(x), \quad 0<x<\infty, \eqno (1.2a)$$
$$ q(0,t) = g_0(t), \quad t>0, \eqno (1.2b)$$
where $q(x,t)$ is a complex-valued function, 
 $u(x)$ is a given real-valued decaying function, and $q_0(x)$,
$g_0(t)$ are given complex-valued functions with sufficient smoothness and
decay which are compatible at $x=t=0$, i.e. $q_0(0) = g_0(0)$.  

(ii) A Neumann problem for the time-dependent \s equation, i.e. the problem
defined by (1.1), (1.2a), and by the equation
$$
q_x(0,t) = g_1(t), \quad t>0, \eqno (1.2b)'
$$
where $\dot{q}_0(0) = g_1(0)$. 

(iii) A Poincar\'e problem for a variable coefficient generalisation of the
Laplace equation on the quarter plane,
$$ q_{xx} + q_{yy} + u(x)q =0, \quad 0<x<\infty, \quad 0<y<\infty, \eqno
(1.3)$$
$$q_y(x,0) + \gamma_1q(x,0) = f(x), \quad 0<x<\infty, \eqno (1.4a)$$
$$ q_x(0,y) + \beta q_y(0,y) + \gamma_2q(0,y) = g(y), \quad 0<y<\infty, \eqno
(1.4b)$$
where $q(x,y)$ is a real-valued function, $u(x)$, $f(x)$, $g(y)$, are given
real-valued functions with appropriate smoothness and decay, and
$\beta,\gamma_1,\gamma_2$ are given real constants.

\paragraph{A. The Case of Constant Coefficients} \ \ 

A new method for studying boundary value problems for linear and for integrable
nonlinear PDEs in two independent variables was introduced in [1].  For linear
PDEs in a convex polygon it involves two novel steps [2]:

(a) {\it Construct an integral representation of $q$ in the complex $k$-plane.}
For example, for a linear dispersive evolution equation on the half line this
representation is [3]
$$ 
q(x,t) = \frac{1}{2\pi} \inta dke^{ikx-i\omega(k)t} \hat q_0(k) + \frac{1}{2\pi}
\int_{\partial D_+} dke^{ikx-i\omega(k)t} \hat g(k,t), \eqno (1.5)
$$
where $\omega(k)$ is the associated dispersion relation, $\hat q_0(k)$ is the
Fourier transform of $q(x,0)$, $\hat g(k,t)$ is an appropriate time
transform of $q(0,t)$ and of $\{\partial_{x}^{l} q(0,t) \}_{1}^{N-1}$
($N$ is the order of the highest spatial derivative), and $\partial
D_+$ is the boundary  of $D_+$ oriented so that $D_+$ is on the left
of $\partial D_+$, where $D_+$ is the domain in the upper half of the
complex $k$-plane defined by $\mathrm{Im}\; \omega(k) \geq 0$.

In the case of equation (1.1) with $u(x) = 0$, it can be shown that
$\omega(k) = ik^2$, $\partial D_+$ is the boundary of the first
quadrant of the complex $k$-plane, and 
$$ 
\hat q_0(k) = \into \! dx e^{-ikx} q_0(x), \quad \mathrm{Im}\; k \leq
0, \eqno (1.6)
$$
$$
\hat g(k,t) = k \hat g_0(k,t) - i\hat g_1(k,t), \quad k \in
\mathbb{C}, \eqno (1.7)
$$
where $\hat{g}_{0}$ and $\hat{g}_{1}$ are defined by
$$ 
\hat g_0(k,t) = \int^t_0 d\tau e^{ik^2\tau} q_0(0,\tau), \quad \hat
g_1(k,t) = \int^t_0 d\tau e^{ik^2\tau} q_x(0,\tau). \eqno (1.8)
$$
We note that in the Dirichlet problem $\hat{g}(k,t)$ involves the {\it unknown}
function $\hat g_1(k,t)$, while in the Neumann problem $\hat{g}(k,t)$ involves the
unknown function $\hat g_0(k,t)$. 

In the case of the Laplace equation on the quarter plane, $q(x,y)$ admits an
integral representation defined along the oriented boundary of the first
quadrant of the complex $k$-plane [2], [4].  For the Dirichlet problem this
representation involves the {\it unknown} Fourier transforms of $q_y(0,x)$ and
of $q_x(0,y)$.

(b) {\it Use the invariant properties of a certain algebraic global relation to
express the unknown transforms in terms of transforms of the given
initial and boundary conditions.}

In the case of equations (1.1), (1.2) with $u(x) =0$, it can be shown
that this yields
$$
-i\hat g_1(k,t) = k\hat g_0(k,t) - \hat q_0(-k) + e^{ik^2t} \hat q(-k,t),
\quad \Im k \geq 0, \eqno (1.9)
$$
where $\hat q(k,t)$ denotes the Fourier transform of $q(x,t)$.  We note that
the rhs of equation (1.9), in addition to the known functions $\hat g_0$  and
$\hat q_0$, it also involves the {\it unknown} function $\hat q(-k,t)$.
However, this term does {\it not} contribute to $q(x,t)$; indeed it gives rise
to the term
$$ 
\int_{\partial D_{+}} dk e^{ikx-ik^2t} \left( e^{ik^2t}\hat
q(-k,t)\right), 
$$
which vanishes according to Jordan's lemma, since both $\exp[ikx] $
and $\hat q(-k,t)$ are analytic and bounded in the first quadrant of
the complex $k$-plane.  Thus dropping the last term of the rhs of
equation (1.9) and substituting $-i\hat{g}_{1}$ in the rhs of equation
(1.7) we find
$$ 
\hat g(k,t) = 2k\hat g_0(k,t) - \hat q_0(-k), \quad \Im k\geq 0,
$$
which involves only transforms of the known functions $q_0(x)$ and $g_0(x)$.

In the case of equations (1.3), (1.4) with 
$u(x)=0$, the analysis of the global relation implies that the unknown transforms can be computed
through the solution of a scalar Riemann--Hilbert problem [4].  For particular values of
$\beta,\gamma_1,\gamma_2$ (which include both the Dirichlet and the Neumann cases) this
Riemann--Hilbert problem can be bypassed and the unknown transforms can be computed using only
algebraic manipulations (just like the case for equations (1.1), (1.2)).

For {\it evolution} PDEs there exist at least three different ways of implementing step (a).
Use: (i) The spectral analysis of the associated Lax pair [2].  (ii) A reformulation of Green's
theorem [5].  (iii) The deformation of the Fourier transform representation [6].  Methods (i) and
(ii) can also be used for {\it elliptic} PDEs.

For particular PDEs, such as equation (1.1) with $u(x) =0$, it is possible to avoid the
construction of an 
integral representation in the complex $k$-plane, and to use instead the usual integral
representations along the Re $\! k$-axis.  However, in general the
complex $k$-plane cannot be avoided; in particular this is the case for an
evolution PDE involving a third order derivative.  

\paragraph{B.  The Case of Variable Coefficients} \ \ 

The analysis of certain classes of PDEs with variable coefficients is also based on the
two novel steps (a) and (b) mentioned in A above.  

\vskip .2in

\noindent (a) {\it The construction of an integral representation.}  For both evolution
and elliptic PDEs this can be achieved using the {\it simultaneous spectral analysis of
the associated Lax pair.}  A Lax pair for equation (1.1) is
$$ \mu_{xx} + (u(x)+k^2) \mu = q \eqno (1.10a)$$
$$ \mu_t + ik^2\mu = iq. \eqno (1.10b)$$
Indeed, equations (1.10) are compatible iff $q(x,t)$ satisfies (1.1): Applying the
operator $\partial_t + ik^2$ to equation (1.10a), noting that  this operator commutes
with $\partial_x^2 + u(x) + k^2$, and using (1.10b) we find (1.1).  Similarly a Lax pair
for equation (1.3) consists of equation (1.10a) together with 
$$ \mu_{yy} - k^2\mu = -q. \eqno (1.11)$$

For evolution PDEs, an alternative approach to constructing integral representations is
to use {\it a completeness relation of the associated space dependent eigenfunctions.}
This is the analogue of the deformation of the Fourier transform, see (iii) in A.  For
equation (1.1) the relevant completeness relation is the restriction to the half line of
the classical relation associated with the time-independent \s equation [7]:  Let $u(x)
\in L_{1}^{2}$, i.e. assume that $\into dx (1+x^2)|u(x)|$ exists. Let $\psi(x,k)$ and
$\phi(x,k)$ be the following solutions of the time-independent \s equation
$$
\chi_{xx} + (u(x) + k^2)\chi = 0, \quad x>0, \quad k \in \mathbb{R}:
\eqno (1.12)
$$
$$ 
\psi(x,k) = e^{ikx} + \frac{1}{k} \int^\infty_x d\xi \sin k (x-\xi) u(\xi)\psi(\xi,k),
\quad x>0, \quad \Im k \geq 0, \eqno (1.13)
$$
$$ 
\phi(x,k) = e^{-ikx} - \frac{1}{k} \int^x_0 d\xi \sin k (x-\xi)u(\xi)\phi(\xi,k),
\quad x >0, \quad k \in \mathbb{C}. \eqno (1.14)
$$
Let $k_j = i p_j$, $j =1,\cdots,n, p_j>0$, be the zeros of $a(k)$,
where $a(k)$ is defined by
$$ 
a(k) = \frac{1}{2ik} (\phi\psi_x - \psi\phi_x), \quad \Im k \geq
0. \eqno (1.15)
$$
Then
$$ 
\delta(x-x') = \frac{1}{2\pi} \inta dk \frac{\psi(x,k)}{a(k)} \phi(x',k) - i \sum^n_1
\frac{\psi(x,ip_j)}{\dot a(ip_j)} \phi(x',ip_j), \eqno (1.16)
$$
where $\dot a(k)$ denotes the derivative of $a(k)$.

\noindent (b) {\it The global relation.}  This equation can be
obtained using the formal adjoint of the given PDE.  For example, if
(1.1) is valid in a simply connected domain with boundary $\Gamma$,
the global relation is
$$
\int_{\Gamma} \Big\{ (qQ) dx + i(q_{x}Q - qQ_{x}) dt \Big\} = 0,
\eqno(1.17)
$$
where $Q$ is any solution of the ``adjoint'' equation
$$
-iQ_{t} + Q_{xx} + uQ = 0. \eqno(1.18)
$$
Indeed, if $q$ and $Q$ satisfy equations (1.1) and (1.18) then
$$
(qQ)_{t} + i(qQ_{x} - q_{x}Q)_{x} = 0,
$$
and Green's theorem yileds (1.17).

For equation (1.13) the global relation is
$$
\int_{\Gamma} \Big\{ (q_{x}Q - qQ_{x}) dy + (qQ_{y} - q_{y}Q)dx \Big\}
  = 0, \eqno(1.19)
$$
where $Q$ is any solution of equation (1.3).

\paragraph{C. Organization of the Paper and Main Results} \ \ 

Using the completeness relation (1.16) as well as the global relation (1.17), it will be
shown in section 2 that the solution of the initial-boundary value problem (1.1), (1.2)
is given by the following formulas:
$$ q(x,t) = \frac{1}{2\pi} \inta dk \frac{\psi(x,k)}{a(k)} \hat q(k,t) - i \sum^n_1
\frac{\psi(x,k_j)}{\dot a(k_j)} \hat q(k_j,t), \eqno (1.20)$$

$$\hat q(k,t) = e^{-ik^2t} \left\{ \hat q_0(k) - \frac{\tilde q_0(k)}{\psi(0,k)} +
\frac{2ka(k)}{\psi(0,k)} \hat g_0(k,t) \right\}, \eqno (1.21)$$
where $\psi(x,k)$, $\phi(x,k)$, $a(k)$ are defined in terms of $u(x)$
by equations (1.13)-(1.15), $\hat g_0(k,t)$ is defined in terms of the
boundary condition $q(0,t) = g_0(t)$ by equation (1.8a), and $\hat
q_0(k)$, $\tilde q_0(k)$ are defined in terms of the initial condition
$q_0(x)$ by
$$ \hat q_0(k) = \into dx \phi(x,k) q_0(x), \quad \tilde q_0(k) = \into dx \psi(x,k)
q_0(x), \quad k \in \mathbb{R}. \eqno (1.22)$$

The solution of the Neumann problem is given by similar formulas but $\hat q(k,t)$ is now
given by
$$
\hat q(k,t) = e^{-ik^2t} \left\{ \hat q_0(k) + ik \frac{\tilde q_0(k)}{\psi_x(0,k)} +
\frac{2a(k)}{\psi_x(0,k)} \hat g_1(k,t) \right\}, \eqno (1.23)
$$
where $\hat{g}_{1}$ is defined in terms of the boundary condition
$q_{x}(0,t) = g_{1}(t)$ by equation (1.8b).

It is possible to deform some of the terms appearing in the integral representation of
$q(x,t)$ to an integral along the boundary of the first quadrant of the complex
$k$-plane. For example it is shown in section 2 that equations (1.20), (1.21) imply
$$ 
q(x,t) = \frac{1}{2\pi} \inta dk \frac{\psi(x,k)}{a(k)} e^{-ik^2t} \hat q_0(k) - i
\sum^n_1 \frac{\psi(x,k_j)}{\dot a(k_j)}  e^{-ik^{2}_{j}t} \hat q_0(k_j)
$$
$$ + \frac{1}{2\pi} \int_{\partial D_1} dk \frac{\psi(x,k)}{\psi(0,k)} e^{-ik^2t} \left[
- \frac{\tilde q_0(k)}{a(k)} + 2k\hat g_0(k,t)\right], \eqno (1.24)$$
where $\partial D_1$ denotes the boundary of the first quadrant of the complex $k$-plane
encircling the points $\{ k_j\}_{1}^{n}$, from the right, see Figure 1.1.

\begin{center}
\begin{minipage}[b]{6cm}
\psfrag{a}{$k_1$}
\psfrag{b}{$k_{n}$}
\psfrag{c}{$\mathrm{Im}\; k$}
\psfrag{d}{$\mathrm{Re}\; k$}
\includegraphics{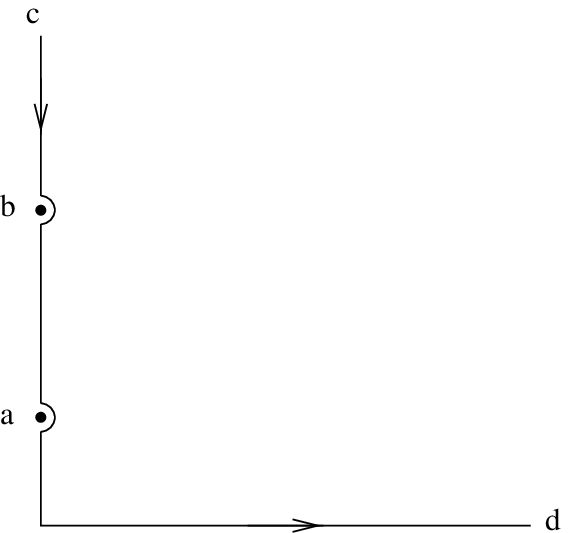}
\centerline{{\bf Figure 1.1:} The contour $\partial{D}_{1}$ associated
  with equations (1.1), (1.2).}
\end{minipage}
\end{center}

\noindent 
If $u(x) = 0$, then $\psi = \exp(ikx)$, $\phi = \exp(-ikx)$, $a(k)=1$, and the
representation (1.24) coincides with the representation (1.5).  The
representation (1.24) is convenient for analyzing the long time
behavior of $q(x,t)$, see section 2.

The particular example of 
$$ u(x) = \frac{2p^2}{\cosh^2p(x-x_0)}, \quad p >0, \quad x_0 >0, \eqno (1.25)$$
is also discussed in section 2; in this case $\psi,\phi,a,$ can be computed explicitly.

In section 3, using the spectral analysis of the Lax pair (1.10a), (1.11), as well as the
global relation (1.19), we will solve equations (1.3), (1.4), see propositions 3.2, 3.3.

In section 4 we discuss further the above results.  In the appendix, for the sake of
completeness, we present a simple derivation of equation (1.16) (see also [7]-[9]).

\section{The Time-Dependent \s Equation}

(a) {\it An Integral Representation} 

The completeness relation (1.16) yields

$$ q(x,t) = \frac{1}{2\pi} \inta dk \frac{\psi(x,k)}{a(k)} \hat q(k,t) - i \sum^n_1
\frac{\psi(x,k_j)}{\dot a(k_j)} \hat q(k_j,t), \eqno (2.1)$$
where
$$\hat q(k,t) = \into dx q(x,t)\phi(x,k),\quad k \in \mathbb{R}. \eqno (2.2)$$
If $q(x,t)$ satisfies equation (1.1) and $\phi(x,k)$ is defined by (1.14) then
$$ \left(e^{ik^2t}q\phi\right)_t + ie^{ik^2t} (q\phi_x-q_x\phi)_x =0. \eqno (2.3)$$
Multiplying equation (2.2) by $e^{ik^2t}$, differentiating with respect to $t$, and using
equation (2.3) we find
$$ \left(e^{ik^2t} \hat q(k,t)\right)_t = ie^{ik^2t} (q(0,t) \phi_x(0,k) - q_x(0,t)
\phi(0,k)).$$
Thus
$$e^{ik^2t} \hat q(k,t) = \hat q(k,0) + i\phi_x(0,k) \hat g_0(k,t) - i\phi(0,k)\hat
g_1(k,t), \eqno (2.4)$$
where $\hat g_0$ and $\hat g_1$ are the $t$-transforms of $q(0,t)$ and
of $q_x(0,t)$ defined by
equations (1.8).

In summary, $q(x,t)$ is given by equation (2.1) where $\hat q(k,t)$ is defined by
equation (2.4).

\vskip .2in
\noindent (b) {\it The Global Relation}

Let $\psi(x,t)$ be defined by equation (1.13), and let
$$ \tilde q(k,t) = \into dxq(x,t)\psi(x,k), \quad \Im k \geq 0. \eqno (2.5)$$
In analogy with equation (2.4) we find
$$ e^{ik^2t} \tilde q(k,t) = \tilde q(k,0) + i\psi_x(0,k) \hat g_0(k,t) -i\psi(0,k) \hat
g_1(k,t), \quad \Im k \geq 0. \eqno (2.6)$$

\paragraph{1. The Dirichlet Case} \ \ 

Solving equation (2.6) for $\hat g_1$ and substituting the resulting expression in (2.4)
we find
$$ \hat q(k,t) = \frac{\phi(0,k)}{\psi(0,k)} \tilde q(k,t) + e^{-ik^2t} \left\{ \hat
q(k,0) - \frac{\phi(0,k)}{\psi(0,k)} \tilde q(k,0) + \frac{2ka(k)}{\psi(0,k)} \hat
g_0(k,t)\right\}. \eqno (2.7)$$
It is important to note that the term $\tilde q(k,t)$ and the associated terms $\tilde
q(k_j,t)$ do {\it not} contribute to $q(x,t)$.  Indeed these terms
give rise to the terms
$$ \frac{1}{2\pi} \inta dk \frac{\psi(x,k)}{a(k)} \frac{\phi(0,k)}{\psi(0,k)} \tilde
q(k,t) - i\sum^n_1 \frac{\psi(x,k)}{\dot a(k)} \frac{\phi(0,k_j)}{\psi(0,k_j)} \tilde
q(k_j,t). \eqno (2.8)$$
The function $\psi(x,k)e^{-ikx}$ is analytic for $\Im k >0$,  thus since $x>0$, the
function $\psi(x,k)$ is also analytic for $\Im k >0$.  Hence all the terms appearing in
the integrant of the first term in (2.8) are analytic for $\Im k>0$.  Furthermore,
$\psi(0,k)\neq 0$ for $\Im k >0$ [8].  Also the contribution from the zeros of $a(k)$ is
cancelled by the second term of equation (2.8).  Hence Jordan's lemma applied in the
upper half of the complex $k$-plane implies that the term (2.8) vanishes.

In summary, $q(x,t)$ is given by (2.1) where $\hat q(k,t)$ is given by equation (2.7)
{\it without} the first term of the rhs of (2.7).

\paragraph{2. The Neumann Case}

Solving equation (2.6) for $\hat g_0$ and proceeding as above we find
$$ \hat q(k,t) = e^{-ik^2t} \left\{ \hat q_0(k) - \frac{\phi_x(0,k)}{\psi_x(0,k)} \tilde
q_0(k) + \frac{2ka(k)}{\psi_x(0,k)} \hat g_1(k^2,t) \right\}. \eqno (2.9)$$

Since $\phi(x,k)$ satisfies (1.14), $\phi(0,k) =1$, $\phi_x(0,k) =-ik$.

\subsection{The Completeness Relation on the Full Line}

There exist a particular class of functions $u(x)$, $-\infty <x<\infty$, called reflectionless
potentials.  For the restriction of this class of potentials on the half line, it is more convenient to
use the completeness relation associated with the full line instead of equation (1.16),
see the Appendix.  This relation is defined by equation (1.16) where $\phi(x,k)$, $a(k)$,
are replaced by $\Phi(x,k)$, $A(k)$, which are defined as follows:
$$
\Phi(x,k) = e^{-ikx} - \frac{1}{k} \int^x_{-\infty} d\xi \sin k(x-\xi) u(\xi)
\Phi(\xi,k), \quad k \in \mathbb{R}, \eqno (2.10)
$$
$$
A(k) = \frac{1}{2ik} (\Phi\psi_x - \Phi_x\psi), \quad \Im k \geq 0. \eqno
(2.11)
$$
The functions $\Phi e^{ikx}$, $A(k)$ are still analytic for $\Im k>0$, furthermore
$\psi(x,k)$ is analytic for $\Im k>0$ provided that $x>0$.  Thus the arguments used for
the derivation of $q(x,t)$ apply mutatis-mutandis if $\phi$, $a$ are replaced by $\Phi,A$.
In summary:

\paragraph{Proposition 2.1}  Let the real-valued function $u(x)$, $-\infty <x<\infty$,
satisfy $\inta dx(1+x^2)|u(x)|<\infty$.  Define $\psi(x,k)$, $\Phi(x,k)$, $A(k)$, in
terms of $u(x)$ by equations (1.13), (2.10), (2.11).  Let $k_j = ip_j$, $p_j>0$,
$j=1,\cdots,n$, be the zeros of $A(k)$ for $\Im k>0$.  The solution of equations (1.1),
(1.2) is given by the following formulas:
$$q(x,t) = \frac{1}{2\pi} \inta dk \frac{\psi(x,k)}{A(k)} \hat q(x,t) - i \sum^n_1
\frac{\psi(x,k_j)}{\dot A(k_j)} \hat q(k_j,t), \eqno (2.12)$$
$$\hat q(k,t) = e^{-ik^2t} \left\{ \hat q_0(k) - \frac{\Phi(0,k)}{\psi(0,k)} \tilde
q_0(k) + \frac{2kA(k)}{\psi(0,k)} \hat g_0(k,t) \right\}, \eqno (2.13)$$
where $\hat g_0$ is defined in terms of the boundary condition $q(0,t)
= g_0(t)$ by equation
(1.8a), while $\tilde q_0(k)$, $\hat q_0(k)$ are defined in terms of the initial
condition $q_0(x)$ by equation (1.22b) and by
$$ \hat q_0(k) = \into dx \Phi(x,k) q_0(x),\quad k \in \mathbb{R}. \eqno (2.14)$$

For the Neumann case
$$ 
\hat q(k,t) = e^{-ik^2t} \left\{ \hat q_0(k) - \frac{\Phi_x(0,k)}{\psi_x(0,k)} \tilde
q_0(k) + \frac{2kA(k)}{\psi_x(0,k)} \hat g_1(k,t)\right\}. \eqno (2.15)
$$

\paragraph{Example 2.1} Let $u(x)$ be given by equation (1.25).  In this case $A(k)$,
$\psi(x,k)$ are defined by equations (A.13), $\Phi(x,k) = a(k)\overline{\psi(x,k)}$, and
$k_1=ip$.  Thus
$$\dot A(k_1) = \frac{1}{2ip}, \quad \psi(x,k_1) = \frac{e^{-px_0}}{2\cosh p(x-x_0)},
\quad \Phi(x,k_1) = \frac{e^{px_0}}{2\cosh p(x-x_0)}.$$
Hence, equations (2.12), (2.13) simplify to the following expressions:
$$ q(x,t) = \frac{1}{2\pi} \left( \inta dk\psi(x,k) \hat Q(k,t) + \psi_1(x) \hat
Q_1(t)\right), \eqno (2.16)$$
$$\hat Q(k,t) = e^{-ik^2t} \left\{ \into dx \overline{\psi(x,k)}q_0(x) - \overline{
\frac{\psi(0,k)}{\psi(0,k)}} \into dx\psi(x,k) q_0(x) + \frac{2k}{\psi(0,k)} \hat
g_0(k,t)\right\}, \eqno (2.17)$$
$$ \psi_1(x) = \frac{\sqrt{p\pi}}{\cosh p(x-x_0)}, \quad \hat Q_1(t) = e^{p^2t} \into dx
\psi_1(x)q_0(x). \eqno (2.18)$$

\subsection{Alternative Representations}

Inserting equation (1.21) into (1.20) we find
\begin{multline}
q(x,t) = \frac{1}{2\pi} \inta dk \frac{\psi(x,k)}{a(k)} e^{-ik^2t}
\hat q_0(k) - i \sum^n_1 \frac{\psi(x,k_j)}{\dot a(k_j)} e^{-ik^2_jt}
\hat q_0(k_j) \\
-\left\{ \frac{1}{2\pi} \inta dk \frac{\psi(x,k)}{\psi(0,k)}
e^{-ik^2t} \frac{\tilde q_0(k)}{a(k)} - i \sum^n_1
\frac{\psi(x,k_j)}{\psi(0,k_j)} e^{-ik^2_jt} \frac{\tilde
  q_0(k_j)}{\dot a(k_j)} \right\} \\
+ \frac{1}{2\pi} \inta dk \frac{\psi(x,k)}{\psi(0,k)}
e^{-ik^{2}t} 2k\hat g_0(k,t). \tag{2.19}
\end{multline}
The term
$$ \frac{\psi(x,k)}{\psi(0,k)} e^{-ik^2t} \frac{\tilde q_0(k)}{a(k)}$$
is bounded and analytic for $\pi/2 < \arg k < \pi$.  Thus using Cauchy's theorem the
integral of this term along the negative real axis can be rotated to an integral along
the positive imaginary axis; the associated residue sum can be incorporated by
encircling the poles from the right.  The integrant of the last integral of the rhs of
(2.19) is also bounded and analytic for $\pi/2 < \arg k < \pi$.  Indeed
$\psi(x,k)/\psi(0,k)$ is analytic for $\Im k>0$ and
$$e^{-ik^2t} \hat g_0(k,t) = \int^t_0 d\tau e^{-ik^2(t-\tau)}g_0(\tau), $$
thus this term is bounded for $k \in[\pi/2,\pi] \cup [3\pi/2,2\pi]$.  Hence the integral
of this term along the negative real axis can also be rotated to an integral along the
positive imaginary axis.  Hence equation (2.19) becomes equation (1.24).

\subsection{The Long $t$ Asymptotics} 

Suppose that equation (1.1) is valid for $0<t<T$. Then $\hat g_0(k,t)$ in (1.24) can be
replaced by $\hat g_0(k,T)$.  Indeed, the difference of these terms gives the
contribution
$$ \frac{1}{2\pi} \int_{\partial D_1} dk \frac{\psi(x,k)}{\psi(0,k)} 2k \int^T_t d\tau
e^{-ik^2(t-\tau)}g_0(\tau),$$
which vanishes due to analyticity considerations.  In particular, if $g_0(t)$ decays as
$t\rightarrow\infty$, $\hat g_0(k,t)$ can be replaced by $\hat g_0(k)$,
$$ \hat g_0(k) = \into d\tau e^{ik^2 \tau}g_0(\tau). \eqno (2.20)$$
Equation (1.24) with $\hat g_0(k,t)$ replaced by $\hat g_0(k)$ has the advantage that
it involves time-dependence only in the explicit exponential form $\exp(-ik^2t)$ and
$\exp(-ik_j^2t)$.  This formulation is most convenient for the analysis of $q(x,t)$ as
$t\rightarrow \infty$, see [10].

\section{The Generalized Laplace Equation}

We first derive an appropriate integral representation for equation (1.3), as well as the
relevant global relations.  We then apply these formulae to the solution of the boundary
value problem (1.3), (1.4).

For the sake of simplicity we assume that the function $a(k)$ has {\it no} zeros for $\Im
k \geq 0$.

\paragraph{Proposition 3.1}  Assume that there exists a real-valued function $q(x,y)$,
$0<x<\infty$, $0<y<\infty$, which satisfies equation (1.3).  Assume that $q(x,y)$ has
sufficient decay as $x\rightarrow  \infty$, $y\rightarrow \infty$, and it also has
sufficient smoothness all the way to the boundary.  Then $q(x,y)$ 
admits the following integral representation:
\begin{multline}
q(x,y) =- \frac{1}{2i\pi} \into dke^{-ky} (\rho(k)\psi(x,k) -
\overline{\rho(k)} \, \overline{\psi(x,k)} ) \\
+ \frac{1}{2i\pi} \int^{i\infty}_0 dk \frac{\psi(x,k)}{\psi(0,k)} (e^{-ky}
g_0(k) + e^{ky}g_0(-k)), \tag{3.1}
\end{multline}
where the functions $g_0(k)$, $\rho(k)$, are defined in terms of the boundary
values $q(x,0)$, $q(0,y)$, by
$$ g_0(k) = \into dye^{ky} q(0,y), \quad \mathrm{Re}\; k \leq 0, \eqno (3.2)$$
$$ \rho(k) = \frac{f_0(k)}{ia(k)}  + \frac{g_0(-k)}{\psi(0,k)}, \quad 0 \leq \arg
k \leq \frac{\pi}{2}, \eqno (3.3)$$
$$ f_0(k) = \into dxq(x,0) \phi(x,k) - \frac{\phi(0,k)}{\psi(0,k)} \into dxq(x,0)
\psi(x,k), \quad \Im k \geq 0. \eqno (3.4)$$
Furthermore, the following global relations are valid
$$ \into dx [kq(x,0) - q_y(x,0)] \phi(x,k) + \into dy e^{ky} [q(0,y) \phi_x(0,k) -
q_x(0,y) \phi(0,k)], \quad \mathrm{Re}\; k \leq 0, \eqno (3.5)$$
$$\into dx [kq(x,0) - q_y(x,0)] \psi(x,k) + \into dy e^{ky} [q(0,y)\psi_x(0,k) -
q_x(0,y)\psi(0,k)], \quad \Im k \geq 0. \eqno (3.6)$$
\paragraph{Proof} \ \

We will construct four functions $\mu_1(x,y,k),\cdots, \mu_4(x,y,k)$ which satisfy
{\it both} equations (1.10a), (1.11), and which for $0<x<\infty$, $0<y<\infty$,
are bounded and analytic in $k$ for $k$ in the first, $\cdots,$ fourth quadrant
of the complex $k$-plane.  Equations (1.10a), (1.11) are invariant under the
transformation $k \rightarrow -k$, thus it is sufficient  to construct only the
functions $\mu_1$ and $\mu_2$, see Figure 3.1.  These functions are given by
\begin{align}
2k\mu_1(x,y,k) &= I_1(x,y,k) + \frac{\psi(x,k)}{\psi(0,k)} I_2(y,k), \quad 0 \leq
\arg k \leq \frac{\pi}{2}, \tag{3.7} \\
2k\mu_2(x,y,k) &= I_1(x,y,k) - \frac{\psi(x,k)}{\psi(0,k)} I_2(y,-k), \quad
\frac{\pi}{2} < \arg k \leq \pi, \tag{3.8}
\end{align}
where
\begin{multline} \tag{3.9} 
I_1(x,y,k) = \frac{1}{ia(k)}\left[ \phi(x,k)\int^\infty_x d\xi
  q(\xi,y) \psi(\xi,k) + \right. \\ 
\left. \psi(x,k)\left( \int^x_0d\xi q(\xi,y)\phi(\xi,k)-
  \frac{\phi(0,k)}{\psi(0,k)}\into d\xi
  q(\xi,y)\psi(\xi,k)\right)\right], \quad \Im k \geq 0, 
\end{multline}
$$I_2(y,k) = \int^y_0 d\eta e^{-k(y-\eta)} q(0,\eta) + \int^\infty_y d\eta
e^{k(y-\eta)} q(0,\eta), \quad {\mathrm{Re}}\; k \geq 0. \eqno (3.10)$$
Indeed, a particular solution of (1.10a) which is analytic for $\Im k>0$, is given
by (see the Appendix)
$$2k\mu = \frac{1}{ia} \left[ \phi\int^\infty_x d\xi q(\xi,y) \psi(\xi,k) + \psi
\int^x_0 d\xi q(\xi,y) \phi(\xi,k)\right]. $$
We add to this particular solution the homogeneous solution $\rho(k,y)\psi(x,k)$
and we choose $\rho(k,y)$ so that at $x=0$ we have an identity:
$$ 2k\mu(x,y,k) = I_1(x,y,k) + 2k \frac{\mu(0,y,k)}{\psi(0,k)} \psi(x,k). \eqno
(3.11)$$
We now use equation (1.11) to determine $\mu(0,y,k)$: the functions
$$ \mu(0,y,k) = \frac{1}{2k} I_2(y,k), \quad \mu(0,y,k) = - \frac{1}{2k}
I_2(y,-k), \eqno (3.12)$$
are particular solutions of equation (1.11) evaluated at $x=0$, and are analytic
for Re $\!k>0$, Re~$\!k<0$, respectively.  Substituting equations (3.12) into (3.11)
we find (3.7), (3.8)

\begin{center}
\begin{minipage}[b]{6cm}
\psfrag{a}{$2k\mu_{2}(x,y,k)$}
\psfrag{b}{$2k\mu_{1}(x,y,k)$}
\psfrag{c}{$2k\mu_{2}(x,y,-k)$}
\psfrag{d}{$2k\mu_{1}(x,y,-k)$}
\includegraphics{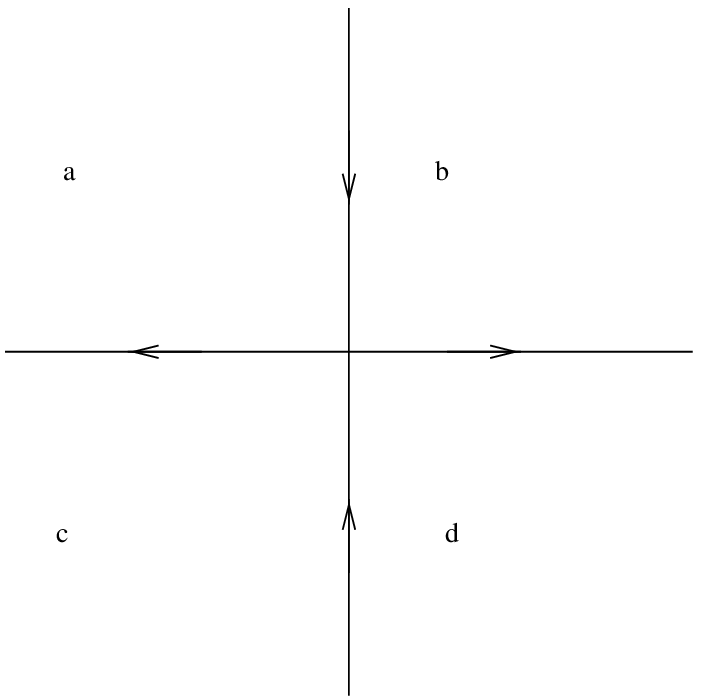}
\centerline{{\bf Figure 3.1:} The domains of analyticity of the
  functions $\mu_1$ and $\mu_2$.}
\end{minipage}
\end{center}

We can use the functions appearing in Figure 3.1 to formulate a Riemann--Hilbert
problem. To this end we need the large $k$ behavior of $\mu_j$, $j=1,2,$ as well
as the relevant ``jumps".  The definitions of $I_1$, $I_2$, imply
\begin{equation}
\begin{split}
k\mu_1(x,y,k) &= O\left(\frac{1}{k}\right), \quad k\rightarrow\infty \quad
{\mathrm{in}} \quad 0< \mathrm{arg} \; k < \frac{\pi}{2}; \\
k\mu_2(x,y,k) &= O\left(\frac{1}{k}\right), \quad k\rightarrow\infty \quad
\mathrm{in} \quad \frac{\pi}{2} < \mathrm{arg}\;  k < \pi. 
\end{split} \tag{3.13}
\end{equation}
In order to compute the jump across the ray $\arg k = \pi/2$ we subtract equations
(3.7), (3.8), and use the definition of $I_2$, 
$$2k\mu_1(x,y,k) - 2k\mu_2(x,y,k) = \frac{\psi(x,k)}{\psi(0,k)} \left[ e^{-ky}
g_0(k) + e^{ky} g_0(-k)\right], \quad \arg k = \frac{\pi}{2}. \eqno (3.14)$$
In order to compute the jump across the ray $\arg k = 0$ we note that since both
$\mu_1(x,y,k)$ and $\mu_2(x,y,-k)$ satisfy equations (1.10b), (1.11), it follows
that their difference satisfies the homogeneous versions of (1.10b),
(1.11).  Also
$\exp[ky]$ cannot appear in this equation since it is unbounded.  Hence
$$2k\mu_1(x,y,k) - 2k\mu_2(x,y,-k) = e^{-ky}
[\rho_1(k)\psi(x,k)+\rho_2(k)\psi(x,-k)], \quad \arg k =0; \eqno (3.15)$$
letting $y=0$, $x\rightarrow \infty$, we find
$$\rho_1(k) = \rho(k), \quad \rho_2(k) = -\overline{\rho(k)}.$$

Let us denote the rhs of equations (3.14), (3.15), by $J_{12}(k)$, $J_{14}(k)$
respectively.  The solution of the above Riemann--Hilbert problem is
$$
2k\mu(x,y,k) = \frac{1}{2i\pi} \left[ - \int^{i\infty}_0dl \frac{J_{12}(l)}{l-k}
+ \into dl \frac{J_{14}(l)}{l-k} - \int^0_{-i\infty} dl \frac{J_{12}(-l)}{l-k} dl -
\int^{-\infty}_0 dl \frac{J_{14}(-l)}{l-k}\right],
$$
or
$$2k\mu(x,y,k) = \frac{1}{2i\pi} \into dlJ_{14}(l) \left[ \frac{1}{l-k} -
\frac{1}{l+k}\right] + \frac{1}{2i\pi} \int^{i\infty}_0 dlJ_{12}(l) \left[
\frac{1}{l+k} - \frac{1}{l-k} \right], \eqno (3.16)$$
where $\mu=\mu_1$ for $0<\arg k < \frac{\pi}{2}$ and $\mu = \mu_2$ for $\pi/2 <
\arg k < \pi$.

Equation (1.10a) implies $q = \lim_{k \rightarrow \infty} (k^2\mu)$; this equation
and equation (3.16) yield
$$ 
q(x,y) =- \frac{1}{2i\pi} \into dl J_{14}(l) + \frac{1}{2i\pi} \int^{\infty}_{0} dl
J_{12} (l), 
$$
which is equation (3.1).

If $\Gamma$ is the boundary of the first quadrant of the complex $z$-plane, the
global relation (1.19) becomes
$$\into dx [q(x,0)Q_y(x,0) - q_y(x,0)Q(x,0)] + \into dy[q(0,y)Q_x(0,y) -
q_x(0,y)Q(0,y)]=0.$$
Letting 
$$ Q=e^{ky}\phi(x,k), \quad Q = e^{ky}\psi(x,k), $$
in the above equation we find equations (3.4), (3.5). \hfill QED

\paragraph{Remark 3.1}  The rhs of equation (3.1) is a real-valued function: The
integrant of the first integral is imaginary, also $\psi(x,k)$ is real for
$k$ purely imaginary thus the integrant of the second integral is
imaginary as well (due to the term $dk$).

The integral representation (3.1) and the global relations (3.4), (3.5) can be
used to solve a variety of boundary value problems.  As an example we consider the
boundary conditions (1.4).

\paragraph{Proposition 3.2} The solution of the boundary value problem (1.3),
(1.4) with $\beta \neq 0$ is given by equation (3.1) where the functions $\rho(k)$
and $g_0(k)$ can be computed in terms of the given functions $f(x)$ and $g(y)$ as
follows: The sectionally analytic function $\{ g_0(k), g_0(-k)\}$ satisfies the
scalar Riemann--Hilbert problem,
$$ g_0(k) \ \ {\mathrm{is \ \ analytic \ \ for \ \ Re}}\; k < 0,$$
$$ g_0(k) = O \left( \frac{1}{k}\right), \quad k \rightarrow \infty,$$
$$ \left\{ \begin{array}{ll} 
g_0(k) + J(k)g_0(-k) = \frac{2k\beta\psi(0,k)}{(k-\gamma_1)D(k)} q(0,0) + F(k), &
\arg k = \frac{\pi}{2}, \\ \\ 
g_0(-k) + J(-k)g_0(k) = \frac{2k\beta\psi(0,-k)}{(k+\gamma_1)D(-k)} q(0,0) +
F(-k), & \arg k = \frac{3\pi}{2}, \end{array} \right. \eqno (3.17)$$
where
$$ 
J(k) = \frac{k+\gamma_1}{k-\gamma_1} \frac{\psi_x(0,k) +
(\gamma_2+\beta k)\psi(0,k)}{\psi_x(0,k) + (\gamma_2-\beta k)\psi(0,k)}, \quad
D(k) = \psi_x(0,k) + (\gamma_2-\beta k)\psi(0,k), \eqno (3.18)
$$
$$ F(k) = \frac{(k-\gamma_1)H(k) + (k+\gamma_1) \overline{H(-\bar k)}}{(k
-\gamma_1)D(k)}, \eqno (3.19)$$
$$H(k) = \into dx f(x)\psi(x,k) + \psi(0,k) \into dye^{ky}g(y).\eqno (3.20)$$
The function $\rho(k)$ is given by
$$\rho(k) =- \frac{k+\gamma_1}{k-\gamma_1} \frac{g_0(-k)}{\psi(0,k)} + N(k), \quad
0 \leq \arg k \leq \frac{\pi}{2}, \eqno (3.21a)$$
$$N(k) =- \frac{\psi(0,k)\into dxf(x)\phi(x,k) - \phi(0,k) \into
dxf(x)\psi(x,k)}{i a(k)(k-\gamma_1)\psi(0,k)}. \eqno (3.21b)$$
\paragraph{Proof} \ \ 

Solving the boundary conditions (1.4) for $q_y(x,0)$, $q_x(0,y)$, substituting the
resulting expressions in (3.5), and integrating by parts the term involving
$q_y(0,y)$, we find
$$ (k+\gamma_1)\into dx q(x,0)\phi(x,k) + d(k)g_0(k) = \beta\phi(0,k)q(0,0) +
h(k), \quad {\mathrm{Re}} \ \ k \leq 0, \eqno (3.22)$$
where the known functions $d(k)$, $h(k)$ are defined by
$$ d(k) = \phi_x(0,k) + (\gamma_2-\beta k)\phi(0,k), \quad h(k) = \into
dxf(x)\phi(x,k) + \phi(0,k)\into dye^{ky}g(y).\eqno (3.23)$$
Similarly, equation (3.6) implies
$$(k+\gamma_1)\into dxq(x,0)\psi(x,k) + D(k)g_0(k) =
\beta\psi(0,k)q(0,0) + H(k),\quad 
\frac{\pi}{2} \leq \arg k \leq \pi, \eqno (3.24)$$
where the known functions $D(k)$, $H(k)$ are defined by (3.18b), (3.20).  Both
equations (3.22), (3.24) are valid for $\pi/2 \leq \arg k \leq \pi$.  Solving these
equations for $f_0(k)$ we find 
$$
f_0(k) = \frac{2ika(k)g_{0}(k)}{\psi(0,k)(k+\gamma_1)} + \frac{\psi(0,k)\into
dxf(x)\psi(x,k) - \phi(0,k)\into dxf(x)\psi(x,k)}{\psi(0,k)(k+\gamma_1)},\;
\frac{\pi}{2} \leq \arg k \leq \pi. 
$$
Taking the complex conjugate of this equation and letting $k\rightarrow -\bar k$
we find an expression for $f_0(k)$ valid for $0\leq  \arg k \leq \pi/2$.
Substituting this expression into equation (3.3) we find (3.21a).

We now show that equation (3.24) implies a Riemann--Hilbert problem for $\{
g_0(k),g_0(-k)\}$.  Taking the complex conjugate of equation (3.24) and letting $k
\rightarrow -\bar k$, we find an equation valid for $0 \leq \arg k \leq \pi/2$, 
$$ (-k+\gamma_1) \into dxq(x,0) \psi(x,k) + \overline{D(-\bar k)} g_0(-k) =
\beta\psi(0,k)q(0,0) + \overline{H(-\bar k)}, 0 \leq \arg k \leq \pi/2. \eqno
(3.25)$$
Eliminating from equations (3.24), (3.25) the integral involving $q(x,0)\psi(x,k)$
we find equation (3.17a).  This equation defines the relevant ``jump" across the
positive imaginary axis; letting $k \rightarrow -k$ we find the jump across the
negative imaginary axis, i.e. equation (3.17b).  The definition of $g_0(k)$
implies $g_0(k) = O(1/k)$, $k\rightarrow \infty$.  \hfill QED

\paragraph{Remark 3.2}  It was shown in [4] that for certain particular cases of
Poincar\'e boundary conditions it is possible to bypass the relevant
Riemann--Hilbert problem.  Such cases were called algebraic.  In what
follows we discuss an algebraic case for equation (1.3).

\paragraph{Proposition 3.3}  The solution of the boundary value problem (1.3),
(1.4) with $\beta =0$, $\gamma_1 <0$, is given by 
\begin{multline} \tag{3.26}
q(x,y) = \frac{1}{2i\pi} \left\{ \int^0_\infty dke^{-ky} \psi(x,k)N(k) +
\int^{i\infty}_0 dke^{-ky} \psi(x,k) \frac{F(k)}{\psi(0,k)} \right\} \\
+ \frac{1}{2i\pi} \left\{ \int^0_{-\infty} dke^{ky} \psi(x,k) \overline{N(-\bar
k)} - \int^{i\infty}_0 dke^{-ky} \psi(x,k) \frac{k-\gamma_1}{k+\gamma_1}
\frac{F(k)}{\psi(0,k)} \right\},
\end{multline}
where $N(k)$ and $F(k)$ are defined in terms of the given functions $f(x)$ and
$g(y)$ by equations (3.21b) and (3.19).

\paragraph{Proof}  The rhs of equation (3.1) can be rewritten in the form
\begin{multline*} \tag{3.27}
q(x,y) = \frac{1}{2i\pi} \left[ \int^0_\infty dke^{-ky} \psi(x,k)\rho(k) +
\int^{i\infty}_0 dke^{-ky} \psi(x,k) \frac{g_0(k)}{\psi(0,k)}\right] \\
+ \frac{1}{2i\pi} \left[ \int^0_{-\infty} dke^{ky} \psi(x,k)
\overline{\rho(-\bar k)} + \int^{i\infty}_0 dke^{ky} \psi(x,k)
\frac{g_0(-k)}{\psi(0,k)} \right].
\end{multline*}
For $\beta=0$, equation (3.17a) becomes
$$ g_0(k) + \frac{k+\gamma_1}{k-\gamma_1} g_0(-k) = F(k), \quad \arg k =
\frac{\pi}{2}. \eqno (3.28)$$
In equation (3.27) we replace $\rho(k)$ by equation (3.21a); also we express
$g_0(k)$ appearing in the first bracket of equation (3.27) in terms of $g_0(-k)$
(using (3.28)), and $g_0(-k)$ appearing in the second bracket of equation (3.27)
in terms of $g_0(k)$ (using again (3.28)).  This yields the rhs of equation
(3.26), and the additional terms
$$ - \frac{1}{2i\pi} \int_{\partial D_1} dke^{-ky} \psi(x,k)
\frac{g_0(-k)}{\psi(0,k)} \frac{k+\gamma_1}{k-\gamma_1} - \frac{1}{2i\pi}
\int_{\partial D_2} dke^{ky} \psi(x,k) \frac{g_0(k)}{\psi(0,k)}
\frac{k-\gamma_1}{k+\gamma_1}, \eqno (3.29)$$
where $\partial D_1$, $\partial D_2$ denote the boundaries of the first and second
quadrants of the complex $k$-plane.  These terms vanish due to the analyticity of
the relevant integrants.  \hfill QED

\paragraph{Remark 3.3}  If $\gamma_1>0$, there exists an additional contribution
due to the pole in the expressions (3.29).  This contribution can be computed
explicitly by evaluating the global relations at $k = -\gamma_1$.  The procedure
is similar to the one discussed in [4].

\section{Discussion and Generalisations}

In this
paper we have introduced a new method for solving boundary value problems for a
class of linear PDEs with variable coefficients, in two independent variables.
This method which can be used to analyze both evolution PDEs (such as equation
(1.1)) as well as elliptic PDEs (such as equation (1.3)), is based on two novel
steps: (a) Construct an integral representation for the solution $q$ in the
complex $k$-plane, see for example equations (1.24) and (3.1).  This
representation involves appropriate transforms of boundary values of $q$ (and of
its derivatives); some of these boundary values are {\it not} prescribed as
boundary conditions.  (b) Use certain algebraic global relations to determine the
transforms of the unknown boundary values; for evolution equations this step
involves only algebra, while for elliptic PDEs, depending on the nature of the
boundary conditions, it involves either algebra (see for example proposition 3.3),
or the solution of a Riemann--Hilbert problem (see for example proposition 3.2). 

The remarkable fact about the new method is that the analysis of the global
relations for the case of PDEs with variable coefficients, is {\it similar} to the
analysis of PDEs with constant coefficients.  This implies that such PDEs can be
solved with the same level of efficiency as the corresponding PDEs of constant
coefficients, except that one now uses appropriate base functions (for example
$\psi(x,k), \phi(x,k)$) instead of the exponential function.  For particular cases
of variable coefficients these base functions can be constructed explicitly. 

For both evolution and elliptic PDEs the relevant integral representations can be
constructed by performing the simultaneous spectral analysis of the associated Lax
pair (see for example Proposition 3.1).  This provides a generalisation to PDEs
with variable coefficients of the formalism introduced in [1].  For evolution PDEs
there exist also an alternative procedure for deriving integral representations;
this is based on the associated completeness relation and on
contour deformation (see for example Section 2).  This approach generalizes to the
case of variable coefficients the formalism introduced in [6]. 

In this paper for economy of presentation we have {\it assumed} that the solution of the
given boundary value problem exists.  However, it is possible to present rigorous
theorems {\it without} the apriori assumption of existence; for the case of
constant coefficients this is discussed in [3] and [11].

We conclude with some remarks. 

1.  Both the Dirichlet and the Neumann problems of equation (1.1) can be solved
using suitable completeness relations instead of the completeness relation (1.16).
For the Dirichlet problem the relevant relation is
$$ 
\delta(x-x') = \frac{1}{2\pi} \inta dk \frac{\psi(x,k)}{a(k)} \left( \phi(x',k)
- \frac{\psi(x',k)}{\psi(0,k)}\right) - i \sum^n_1 \frac{\psi(x,k_j)}{\dot
a(k_j)} \left( \phi(x',k_j) - \frac{\psi(x',k_j)}{\psi(0,k_j)}\right). \eqno
(4.1)
$$
This relation can be derived by performing the spectral analysis of equation
(1.1a) with the condition $\mu(0,k) =0$ [9].  Using (4.1) it follows that instead
of equation (2.4) we now find
$$e^{ik^2t} \hat q(k,t) = \hat q(k,0) + k\hat g_0(k,t).$$
Thus we find $\hat q(k,t)$ immediately without the need to use the global
relation.  Equation (4.1) is the generalisation of the sine-transform, while equation (1.16)
is the generalisation of the Fourier transform (restricted on the half-line).  It
has been emphasized by the author that there does {\it not} exist an analogue of
the sine transform for problems with derivatives of third order.  For example the
spectral analysis of the equation
$$ 
\frac{d^3f}{dx^3} + \lambda f =0, \quad 0<x<\infty, \quad
\mathrm{with} \; f(0,\lambda) =0,
$$
in Rietz is incomplete [12].  An advantage of the new method introduced by the author
is that it does {\it not} rely on the spectral analysis of the associated
differential operator obtained by separation of variables.  Similar considerations
apply to PDEs with variable coefficients.  Consider for example the Dirichlet
problem for the PDE
$$
q_t + q_{xxx} + u(x)q=0, \quad 0<x<\infty. \eqno (4.2)
$$
For this equation there does {\it not} exist an analogue of equation
(4.1), but using the spectral theory developed in [13], [14], it
\emph{is} possible to obtain the analogue of equation (1.16).  This
implies that equation (4.2) {\it can} be
solved using the method introduced in this paper, but {\it cannot} be solved by
performing the spectral analysis of
$$\mu_{xxx} + (u(x) + k^3)\mu = q(x), \quad 0<x<\infty,$$
and demanding $\mu(0,k) =0.$

2. We note that the classical method of separation of variables {\it fails} for
the boundary value problem (1.3), (1.4).  Indeed, although both the PDE and the
domain are separable, the boundary conditions at $x=0$ are {\it not} separable.
For example, in the particular case of $\gamma_1=0$ the boundary condition at
$y=0$ is consistent with the cosine transform,
$$ \hat q(x,k) = \into dyq(x,y)\cos ky.$$
The cosine transform of equation (1.3) yields
$$ \hat q_{xx} + u(x) \hat q - k^2\hat q = f(x).$$
But the cosine transform of the boundary condition at $x=0$ leads to the {\it
non-separable} equation
$$ \hat q_x(0,k) + \gamma_2\hat q (0,k) + \beta \into dyq_y(0,y) \cos ky = \into
dyg(y)\cos ky.$$

3. Using the method introduced in this paper it is possible to analyze a large
class of other boundary value problems.  As an example we mention the following
generalisation of the Helmholtz equation on the quarter plane,
$$ q_{xx} + q_{yy} + \alpha q + (u(x) + v(y))q=0, \quad 0<x<\infty, \quad
0<y<\infty, \eqno (4.3)$$
where $\alpha$ is a real constant and $u(x),v(y),$ are given real-valued functions
with sufficient decay as $x\rightarrow \infty$, $y\rightarrow \infty$.

4. The implementation of the new method to linear evolution PDEs with
constant coefficients in {\it two space dimensions} is carried out in [6].
In this case one constructs an integral representation of $q(x_1,x_2,t)$ in
the complex $(k_1,k_2)$-planes instead of the complex $k$-plane.  Using
the results of [6] together with the methodology introduced  in \S 2, it should
be possible to solve boundary value problems for certain classes of evolution PDEs in
two space dimensions with variable coefficients.  The main difference
is that, in order to derive an  integral representation for
$q(x_1,x_2,t)$, one must now use an appropriate completeness relation
instead of the two-dimensional Fourier transform.  In recent years
progress has been made for constructing such completeness relations
for two dimensional eigenvalue equations, see for example [15], [16].

\section*{Appendix -- Completeness Relations}
Given $u(x)$, define $\psi(x,k)$, $\phi(x,k)$, $a(k)$, by equations (1.13)-(1.15).
Evaluating the rhs of equation (1.15) at $x=0$, and as $x\rightarrow \infty$, it
follows that
$$ 
a(k) = 1 + \frac{1}{2ik} \into d\xi e^{ik\xi} u(\xi) \phi(\xi,k) =
\frac{1}{2ik} \Big(\psi_x(0,k) + ik\psi(0,k) \Big). \eqno (A.1)
$$
The definition of $\psi(x,k)$ implies 
$$ 
\psi(x,k) e^{-ikx} = 1 + \frac{1}{2ik} \int^\infty_x d\xi \left(
1-e^{2ik(\xi-x)}\right) u(\xi) \psi(\xi,k) e^{-ik\xi}, \eqno (A.2)
$$
thus $\psi(x,k)e^{-ikx} $ is analytic for $\Im k >0$.  The function $\phi(x,k)$ is
an entire function of $k$; writing an equation analogous to equation (A.2) for the
function $\phi(x,k)e^{ikx}$ it follows that this function is bounded for $\Im k
\geq  0 $.  Equations (A.1) imply that $a(k)$ is also analytic for $\Im k >0$.  In
summary 
$$ \psi(x,k)e^{-ikx}, \quad \phi(x,k)e^{ikx}, \quad a(k) \eqno (A.3)$$
are bounded and analytic for $\Im k>0$. 

We now consider equation
$$\mu_{xx} + (u(x)+k^2)\mu = q(x), \quad 0<x<\infty, \quad k \in \mathbb{R}, \eqno
(A.4)$$
where $q(x)$ has sufficient smoothness as well as sufficient decay as
$x\rightarrow \infty$.  Using variation of parameters it follows that for $k\in
\mathbb{R}$, a particular solution of equation (A.4) is
$$ 2ika(k)\mu_p = -\phi(x,k) \int^xd\xi q(\xi)\psi(\xi,k) + \psi(x,k) \int^xd\xi
q(\xi)\phi(\xi,k).$$
We choose the limits of integration in such a way that the rhs of the above
equation is analytic for $\Im k >0$: 
$$ 2ika(k)\mu^+(x,k) = \phi(x,k) \int^\infty_x d\xi q(\xi)\psi(\xi,k) + \psi(x,k)
\int^x_0d\xi q(\xi)\phi(\xi,k), \quad \Im k \geq 0. \eqno (A.5)$$
Indeed, the first term of the rhs of (A.5) equals
$$ 
(\phi(x,\xi)e^{ikx}) \int^\infty_x d\xi
e^{ik(\xi-x)}q(\xi)(\psi(\xi,k)e^{-ik\xi}),
$$
and $\exp[ik(\xi-x)]$ as well as the two terms appearing in the parentheses are
analytic for $\Im k>0$, see (A.3); similarly for the second term of the rhs of
equation (A.5).  

Both functions $\mu^+(x,k)$ and $\mu^+(x,-k)$ satisfy the same equation (A.4),
thus their difference satisfies the homogeneous version of equation (A.4), hence
$$2ik\mu^+(x,k) - 2ik\mu^+(x,-k) = \rho(k)\psi(x,k) + \rho(-k)\psi(x,-k), \quad k
\in \mathbb{R}. \eqno (A.6a)$$
Evaluating this equation as $x\rightarrow \infty$, it follows that 
$$ 
\rho(k) = \frac{1}{a(k)}\into d\xi q(\xi) \phi(\xi,k). \eqno (A.6b)
$$
Equation (A.5) implies that $\mu^+(x,k)$ is a meromorphic function of $k$ for $\Im
k >0$.  Also using the definitions of $\phi$ and $\psi$ it follows that
$2ik\mu^+(x,k) = O(1/k)$ as $k\rightarrow \infty$.  These facts together with
equation (A.6a) define a Riemann--Hilbert problem for the sectionally meromorphic
function $\{ 2ik\mu^+(x,k),2ik\mu^+(x,-k)\}$.  Before solving this problem we need
to compute the relevant residues.  Let $k_j = ip_j$, $p_j>0$, $j=1,\cdots,n$, be
the zeros of $a(k)$ for $\Im k>0$; the residue of $2ik\mu^+(x,k)$ at $k=k_j$ is 
$$\frac{1}{\dot a(k_j)} \left\{ \phi(x,k_j)\int^\infty_x d\xi q(\xi)\psi(\xi,k_j)
+ \psi(\xi,k_j) \int^x_0 d\xi q(\xi)\phi(\xi,k_j)\right\} = \frac{\psi(x,k_j)
\into d\xi q(\xi)\phi(\xi,k_j)}{\dot  a(k_j)}, \eqno (A.7)$$
where we have used $\phi(x,k_j) = c_j\psi(x,k_j)$ (which follows from
equation (1.15), since $a(k_j) =0$).  Equation (A.6a) can be rewritten as
\begin{multline*}
\left\{ 2ik\mu^+(x,k) - \sum^n_1 \frac{A_j(x)}{k-k_j}\right\} - \left\{
2ik\mu^+(x,-k) - \sum^n_1 \frac{A_j(x)}{k+k_j} \right\} = \\
\rho(k)\psi(x,k) + \rho(-k)\psi(x,-k) - \sum^n_1 \frac{A_j(x)}{k-k_j} +
\sum^n_1 \frac{A_j(x)}{k+k_j},
\end{multline*}
where $A_j(x)$ denotes the rhs of equation (A.7). This equation implies
$$
2ik\mu^+(x,k) - \sum^n_1 \frac{A_j(x)}{k-k_j} = \frac{1}{2i\pi} \inta dl
\frac{\rho(l)\psi(x,l) + \rho(-l)\psi(x,-l)}{l-k} + \sum^n_1
\frac{A_j(x)}{k+k_j}, \quad \Im k >0. \eqno (A.8)
$$
Using $q = \lim_{k\rightarrow \infty} (k^2\mu)$, equation (A.8) yields 
$$ q(x) = \frac{1}{2\pi} \inta dk \psi(x,k) \rho(k) - i\sum^n_1 A_j(x),$$
which is equivalent to equation (1.16).

\paragraph{Remark A.1}  $\phi(x,k), \psi(x,k)$ are related by
$$ 
\phi(x,k) = a(k)\psi(x,-k) + b(k)\psi(x,k), \quad k \in \mathbb{R}, \eqno
(A.9)
$$
where
$$ 2ikb(k) = \psi(x,-k)\phi_x(x,k) - \phi(x,k) \psi_{x}(x,-k), \quad k \in \mathbb{R}.
\eqno (A.10)$$
The functions $a(k), b(k)$ satisfy
$$a(k)a(-k) - b(k)b(-k) = 1. \eqno (A.11)$$
Since $u(x)$ is real,
$$ \psi(x,-k) = \overline{\psi(x,\bar k)}, \quad a(-k) = \overline{a(\bar k)},
\quad b(-k) = \overline{b(\bar k)}. \eqno (A.12)$$

\paragraph{Remark A.2}  It is sometimes convenient to consider extensions of
$u(x)$ from $0<x<\infty$ to $-\infty <x<\infty$, see the example below.  The
completeness relation for $x$ on the full line is also given by equation (1.16),
but with $\phi(x,k)$, $a(k)$ replaced by $\Phi(x,k)$, $A(k)$, which we defined by
equations (2.10), (2.11).

\paragraph{Example}  Let $u(x), -\infty <x<\infty$, be defined by equation (1.25).
It is well known that in this case
$$ 
A(k) = \frac{k-ip}{k+ip}, \quad B(k) =0, \quad \psi(x,k) = \frac{k+ip \tanh
p(x-x_0)}{k+ip} e^{ikx}, \eqno (A.13)
$$
where $B(k)$ is defined by equation (A.10) with $\phi$ replaced by
$\Phi$.  On the other hand, if we restrict $u(x)$ to the half line, it follows that
$\psi(x,k)$ is still given by (A.13c) but $a(k), b(k)$, are given by the more
complicated formulae
$$ a(k) = \frac{\psi_x(0,k) + ik\psi(0,k)}{2ik}, \quad b(k) =- \frac{\psi_x(0,-k)
+ ik\psi(0,-k)}{2ik}. $$
Thus if $u(x)$ is given by (1.25) it is more convenient to use the completeness
relation of the full line.

\section*{Acknowledgment}

I am grateful to C. Rogers and W. Schiff for important suggestions.  This
work was supported by the EPSRC.

\end{document}